\documentclass[twoside,12pt,leqno]{article}
\advance\topmargin by -\headheight\advance\topmargin by -\headsep

\usepackage{amsmath}
\usepackage{amsthm}
\usepackage{amssymb}
\usepackage{amscd}
\usepackage{graphicx}
\usepackage{epsfig}

\numberwithin{equation}{section}
\allowdisplaybreaks

\newtheorem{theorem}{Theorem}[section]

\theoremstyle{definition}
\newtheorem{definition}{Definition}[section]
\newtheorem{example}{Example} [section]

\newtheorem{remark}{Remark}[section]
\newtheorem*{remarku}{Remark}
\begin{document}
%%%%%%%%%%%%%%%%%%%%%%%%%%%%%%%%%%%%%%%%%%
%%%%%%%%%%%%%%%%%%%%%%%%%%%%%%%%%%%%%%%%%%%%%%%%
%%%%%%%%  START OF THE CONTRIBUTION
%%%%%%%%%%%%%%%%%%%%%%%%%%%%%%%%%%%%%%%%%%%%
%%%%    private macros
%%%%    do not use renewcommand (i.e. do not redefine
%%%%%   standard  Latex-Commands
%%%%%%%%%%%%%%%%%%%%%%%%%%%%%%%%%%%%%%%%%%%%%%%%%%%%%%%
\font\black=cmbx10 \font\sblack=cmbx7 \font\ssblack=cmbx5
\font\blackital=cmmib10  \skewchar\blackital='177
\font\sblackital=cmmib7 \skewchar\sblackital='177
\font\ssblackital=cmmib5 \skewchar\ssblackital='177
\font\sanss=cmss12 \font\ssanss=cmss8 %scaled 900
\font\sssanss=cmss8 scaled 600 \font\blackboard=msbm10
\font\sblackboard=msbm7 \font\ssblackboard=msbm5
\font\caligr=eusm10 \font\scaligr=eusm7 \font\sscaligr=eusm5
\font\blackcal=eusb10 \font\fraktur=eufm10 \font\sfraktur=eufm7
\font\ssfraktur=eufm5 \font\blackfrak=eufb10

\font\bsymb=cmsy10 scaled\magstep2
\def\all#1{\setbox0=\hbox{\lower1.5pt\hbox{\bsymb
       \char"38}}\setbox1=\hbox{$_{#1}$} \box0\lower2pt\box1\;}
\def\exi#1{\setbox0=\hbox{\lower1.5pt\hbox{\bsymb \char"39}}
       \setbox1=\hbox{$_{#1}$} \box0\lower2pt\box1\;}

\def\mi#1{{\fam1\relax#1}}
\def\tx#1{{\fam0\relax#1}}

\newfam\bifam
\textfont\bifam=\blackital \scriptfont\bifam=\sblackital
\scriptscriptfont\bifam=\ssblackital
\def\bi#1{{\fam\bifam\relax#1}}

\newfam\blfam
\textfont\blfam=\black \scriptfont\blfam=\sblack
\scriptscriptfont\blfam=\ssblack
\def\rbl#1{{\fam\blfam\relax#1}}

\newfam\bbfam
\textfont\bbfam=\blackboard \scriptfont\bbfam=\sblackboard
\scriptscriptfont\bbfam=\ssblackboard
\def\bb#1{{\fam\bbfam\relax#1}}

\newfam\ssfam
\textfont\ssfam=\sanss \scriptfont\ssfam=\ssanss
\scriptscriptfont\ssfam=\sssanss
\def\sss#1{{\fam\ssfam\relax#1}}

\newfam\clfam
\textfont\clfam=\caligr \scriptfont\clfam=\scaligr
\scriptscriptfont\clfam=\sscaligr
\def\cl#1{{\fam\clfam\relax#1}}

\newfam\frfam
\textfont\frfam=\fraktur \scriptfont\frfam=\sfraktur
\scriptscriptfont\frfam=\ssfraktur
\def\fr#1{{\fam\frfam\relax#1}}

\def\cb#1{\hbox{$\fam\gpfam\relax#1\textfont\gpfam=\blackcal$}}

\def\hpb#1{\setbox0=\hbox{${#1}$}
    \copy0 \kern-\wd0 \kern.2pt \box0}
\def\vpb#1{\setbox0=\hbox{${#1}$}
    \copy0 \kern-\wd0 \raise.08pt \box0}

\def\pmb#1{\setbox0\hbox{${#1}$} \copy0 \kern-\wd0 \kern.2pt \box0}
\def\pmbb#1{\setbox0\hbox{${#1}$} \copy0 \kern-\wd0
      \kern.2pt \copy0 \kern-\wd0 \kern.2pt \box0}
\def\pmbbb#1{\setbox0\hbox{${#1}$} \copy0 \kern-\wd0
      \kern.2pt \copy0 \kern-\wd0 \kern.2pt
    \copy0 \kern-\wd0 \kern.2pt \box0}
\def\pmxb#1{\setbox0\hbox{${#1}$} \copy0 \kern-\wd0
      \kern.2pt \copy0 \kern-\wd0 \kern.2pt
      \copy0 \kern-\wd0 \kern.2pt \copy0 \kern-\wd0 \kern.2pt \box0}
\def\pmxbb#1{\setbox0\hbox{${#1}$} \copy0 \kern-\wd0 \kern.2pt
      \copy0 \kern-\wd0 \kern.2pt
      \copy0 \kern-\wd0 \kern.2pt \copy0 \kern-\wd0 \kern.2pt
      \copy0 \kern-\wd0 \kern.2pt \box0}

\def\cdotss{\mathinner{\cdotp\cdotp\cdotp\cdotp\cdotp\cdotp\cdotp
        \cdotp\cdotp\cdotp\cdotp\cdotp\cdotp\cdotp\cdotp\cdotp\cdotp
        \cdotp\cdotp\cdotp\cdotp\cdotp\cdotp\cdotp\cdotp\cdotp\cdotp
        \cdotp\cdotp\cdotp\cdotp\cdotp\cdotp\cdotp\cdotp\cdotp\cdotp}}

%%%%%%%%%%%%%%%%%%%%%%%%%%%%%%%%%%%%%%%%%%%%%%%%%%%%%%%%%%%%
%\textwidth15.6cm \textheight24cm \hoffset-1.6cm \voffset-2.5cm
\font\frak=eufm10 scaled\magstep1 \font\fak=eufm10 scaled\magstep2
\font\fk=eufm10 scaled\magstep3 \font\scriptfrak=eufm10
\font\tenfrak=eufm10

%\newtheorem{theorem}{Theorem}
%\newtheorem{corollary}{Corollary}
%\newtheorem{proposition}{Proposition}
%\newtheorem{definition}{Definition}
%\newtheorem{lemma}{Lemma}
%\font\frak=eufm10 scaled\magstep1
%%%%%%%%%%%%%%%%%%%%%%%%%%%%%%%%%%%%%%%%%%%%%%%%%%%%
%\newenvironment{pf}{{\noindent{\it Proof. }}}{\ \rule{2mm}{2.5mm}\medskip}
%\newenvironment{pft}{{\noindent{\it Proof of Theorem }}}{\ \rule{2mm}
%{2.5mm}\medskip}
%%%%%%%%%%%%%%%%%%%%%%%%%%%%%%%%%%%%%%%%%%%%%%%%%%%%%%%%%%%%%%%%%

\mathchardef\za="710B  %\alpha
\mathchardef\zb="710C  %\beta
\mathchardef\zg="710D  %\gamma
\mathchardef\zd="710E  %\delta
\mathchardef\zve="710F %\epsilon
\mathchardef\zz="7110  %\zeta
\mathchardef\zh="7111  %\eta
\mathchardef\zvy="7112 %\theta
\mathchardef\zi="7113  %\iota
\mathchardef\zk="7114  %\kappa
\mathchardef\zl="7115  %\lambda
\mathchardef\zm="7116  %\mu
\mathchardef\zn="7117  %\nu
\mathchardef\zx="7118  %\xi
\mathchardef\zp="7119  %\pi
\mathchardef\zr="711A  %\rho
\mathchardef\zs="711B  %\sigma
\mathchardef\zt="711C  %\tau
\mathchardef\zu="711D  %\upsilon
\mathchardef\zvf="711E %\phi
\mathchardef\zq="711F  %\chi
\mathchardef\zc="7120  %\psi
\mathchardef\zw="7121  %\omega
\mathchardef\ze="7122  %\varepsilon
\mathchardef\zy="7123  %\vartheta
\mathchardef\zf="7124  %\varomega
\mathchardef\zvr="7125 %\varrho
\mathchardef\zvs="7126 %\varsigma
\mathchardef\zf="7127  %\varphi
\mathchardef\zG="7000  %\Gamma
\mathchardef\zD="7001  %\Delta
\mathchardef\zY="7002  %\Theta
\mathchardef\zL="7003  %\Lambda
\mathchardef\zX="7004  %\Xi
\mathchardef\zP="7005  %\Pi
\mathchardef\zS="7006  %\Sigma
\mathchardef\zU="7007  %\Upsilon
\mathchardef\zF="7008  %\Phi
\mathchardef\zW="700A  %\Omega

\newcommand{\be}{\begin{equation}}
\newcommand{\ee}{\end{equation}}
\newcommand{\ra}{\rightarrow}
\newcommand{\lra}{\longrightarrow}
\newcommand{\bea}{\begin{eqnarray}}
\newcommand{\eea}{\end{eqnarray}}
\newcommand{\beas}{\begin{eqnarray*}}
\newcommand{\eeas}{\end{eqnarray*}}
\newcommand{\R}{{\mathbb R}}
\newcommand{\T}{{\mathbb T}}
\newcommand{\C}{{\mathbb C}}
\newcommand{\1}{{\mathbf 1}}
\newcommand{\SL}{SL(2,\C)}
\newcommand{\Sl}{sl(2,\C)}
\newcommand{\SU}{SU(2)}
\newcommand{\su}{su(2)}
\newcommand{\G}{{\goth g}}
\newcommand{\D}{{\rm d}}
\newcommand{\Df}{{\rm d}^\zF}
\newcommand{\de}{\,{\stackrel{\rm def}{=}}\,}
\newcommand{\we}{\wedge}
\newcommand{\nn}{\nonumber}
\newcommand{\ot}{\otimes}
\newcommand{\s}{{\textstyle *}}
\newcommand{\ts}{T^\s}
\newcommand{\oX}{\stackrel{o}{X}}
\newcommand{\oD}{\stackrel{o}{D}}
\newcommand{\obD}{\stackrel{o}{\bD}}
%%%%%%%%%%%%%%%%%%%%%%%%%%%%%%%%%%%%%%%%%%%%%%%%%%%%%%%%%%%%%%%%%%
\newcommand{\pa}{\partial}
\newcommand{\ti}{\times}
\newcommand{\A}{{\cal A}}
\newcommand{\Li}{{\cal L}}
\newcommand{\ka}{\mathbb{K}}
\newcommand{\find}{\mid}
\newcommand{\ad}{{\rm ad}}
\newcommand{\rS}{]^{SN}}
\newcommand{\rb}{\}_P}
\newcommand{\p}{{\sf P}}
\newcommand{\h}{{\sf H}}
\newcommand{\X}{{\cal X}}
\newcommand{\I}{\,{\rm i}\,}
\newcommand{\rB}{]_P}
\newcommand{\Ll}{{\pounds}}
\def\lna{\lbrack\! \lbrack}
\def\rna{\rbrack\! \rbrack}
\def\rnaf{\rbrack\! \rbrack_\zF}
\def\rnah{\rbrack\! \rbrack\,\hat{}}
\def\lbo{{\lbrack\!\!\lbrack}}
\def\rbo{{\rbrack\!\!\rbrack}}
\def\lan{\langle}
\def\ran{\rangle}
\def\zT{{\cal T}}
\def\tU{\tilde U}
\def\ati{{\stackrel{a}{\times}}}
\def\sti{{\stackrel{sv}{\times}}}
\def\aot{{\stackrel{a}{\ot}}}
\def\sati{{\stackrel{sa}{\times}}}
\def\saop{{\stackrel{sa}{\op}}}
\def\bwa{{\stackrel{a}{\bigwedge}}}
\def\svop{{\stackrel{sv}{\oplus}}}
\def\saot{{\stackrel{sa}{\otimes}}}
\def\cti{{\stackrel{cv}{\times}}}
\def\cop{{\stackrel{cv}{\oplus}}}
\def\dra{{\stackrel{\xd}{\ra}}}
\def\bdra{{\stackrel{\bd}{\ra}}}
\def\bAff{\mathbf{Aff}}
\def\Aff{\sss{Aff}}
\def\bHom{\mathbf{Hom}}
\def\Hom{\sss{Hom}}
\def\bt{{\boxtimes}}
\def\sot{{\stackrel{sa}{\ot}}}
\def\bp{{\boxplus}}
\def\op{\oplus}
\def\bwak{{\stackrel{a}{\bigwedge}\!{}^k}}
\def\aop{{\stackrel{a}{\oplus}}}
\def\V{{\cal V}}
\def\cD{{\cal D}}
\def\cL{{\cal L}}
\def\cR{{\cal R}}
\def\cJ{{\cal J}}
\def\bA{\mathbf{A}}
\def\bI{\mathbf{I}}
\def\wh{\widehat}
\def\wt{\widetilde}
\def\ol{\overline}
\def\Sec{\sss{Sec}}
\def\Lin{\sss{Lin}}
\def\ader{\sss{ADer}}
\def\ado{\sss{ADO^1}}
\def\adoo{\sss{ADO^0}}
\def\AS{\sss{AS}}
\def\bAS{\sss{AS}}
\def\bLS{\sss{LS}}
\def\bAP{\sss{AV}}
\def\bLP{\sss{LP}}
\def\AP{\sss{AP}}
\def\LP{\sss{LP}}
\def\LS{\sss{LS}}
\def\Z{\mathbf{Z}}
\def\oZ{\overline{\bZ}}
\def\oA{\overline{\bA}}
\def\cim{{C^\infty(M)}}
\def\de{{\cal D}^1}
\def\la{\langle}
\def\ran{\rangle}
%%%%%%%%%%%%%%%%%%%
\def\by{{\bi y}}
\def\bs{{\bi s}}
\def\bc{{\bi c}}
\def\bd{{\bi d}}
\def\bh{{\bi h}}
\def\bD{{\bi D}}
\def\bY{{\bi Y}}
\def\bX{{\bi X}}
\def\bL{{\bi L}}
\def\bV{{\bi V}}
\def\bW{{\bi W}}
\def\bS{{\bi S}}
\def\bT{{\bi T}}
\def\bC{{\bi C}}
\def\bE{{\bi E}}
\def\bF{{\bi F}}
\def\bP{{\bi P}}
\def\bp{{\bi p}}
\def\bz{{\bi z}}
\def\bZ{{\bi Z}}
\def\bq{{\bi q}}
\def\bQ{{\bi Q}}
\def\bx{{\bi x}}

\def\sA{{\sss A}}
\def\sC{{\sss C}}
\def\sD{{\sss D}}
\def\sG{{\sss G}}
\def\sH{{\sss H}}
\def\sI{{\sss I}}
\def\sJ{{\sss J}}
\def\sK{{\sss K}}
\def\sL{{\sss L}}
\def\sO{{\sss O}}
\def\sP{{\sss P}}
\def\sPh{{\sss P\sss h}}
\def\sT{{\sss T}}
\def\sV{{\sss V}}
\def\sR{{\sss R}}
\def\sS{{\sss S}}
\def\sE{{\sss E}}
\def\sF{{\sss F}}
\def\st{{\sss t}}
\def\sg{{\sss g}}
\def\sx{{\sss x}}
\def\sv{{\sss v}}
\def\sw{{\sss w}}
\def\sQ{{\sss Q}}
\def\sj{{\sss j}}
\def\sq{{\sss q}}
\def\xa{\tx{a}}
\def\xc{\tx{c}}
\def\xd{\tx{d}}
\def\xi{\tx{i}}
\def\xD{\tx{D}}
\def\xV{\tx{V}}
\def\xF{\tx{F}}

%%%%%%%%%%%%%%%%%%%%%%%%
%%%%%  end of private macros

\setcounter{page}{1}
\thispagestyle{empty}
%\FirstPageHead{14}{2003}{\pageref{firstpage}--\pageref{lastpage}}
%  Parameters: Volume, year, page range,

%%%%%%%%%%%%%%%%%%%%%%%%%%%%%%%%%%%%%%%%%%%%%%%%%%%%%%
%%%  Here  the top-matter of your Article starts
%%%%%%%%%%%%%%%%%%%%%%%%%%%%%%%%%%%%%%%%%%%%%%%%%%%%%%

%%%%    Please replace by your data

%%%   This are the running heads
\markboth{K.~Grabowska, J.~Grabowski,
P.~Urba\'nski}{Frame-independent mechanics}

\label{firstpage}
$ $
\bigskip

\bigskip

%%%%    Title
\centerline{{\Large \bf Frame-independent mechanics: }}

\medskip\noindent
\centerline{{\Large \bf geometry on affine
bundles}}\footnote{Research supported by the Polish Ministry of
Scientific Research and Information Techno\-logy under the grant
No. 2 P03A 036 25.}

%%%   Author
\bigskip
\bigskip
\centerline{{\large Katarzyna Grabowska, Janusz Grabowski, and
Pawe\l \ Urba\'nski}}

\vspace*{.7cm}

\begin{abstract}

Main ideas of the differential geometry on affine bundles are
presented. Affine counterparts of Lie algebroid and Poisson
structures are introduced and discussed. The developed concepts
are applied in a frame-independent formulation of the
time-dependent and the Newtonian mechanics.
\end{abstract}
\begin{quote}

\vspace{-5pt}

\end{quote}

\pagestyle{myheadings}
%%%%%%%%%%%%%%%%%%%%%%%%%%%%%%%%%
\section{Introduction}

    Frame-independent formulation of the Newtonian analytical mechanics can be achieved,
    like in the case of a relativistic charged particle, by increasing the dimension of the configuration space
    (cf. the Kaluza theory).
    There is an alternative approach, based on ideas of Tulczyjew, in which the four-dimensional
    space-time is used as the configuration space. The phase space is no longer a cotangent bundle,
    but an affine bundle, modeled on the cotangent bundle. Also time-dependent mechanics requires affine
    objects. Here Lagrangian is a function on a space of first-jets of motions, interpreted as an affine
    subbundle of the tangent bundle to the space-time. The phase bundle is a vector bundle, however the
    Hamiltonian is not a function, but a section of a bundle over the phase manifold. The structure
    on the phase bundle, which makes possible the Hamiltonian formulation of the dynamics, is no longer
    a standard Poisson structure, but its affine counterpart.

    It is known that Lagrangian formulation of the dynamics of an autonomous system is based on the
    Lie algebroid structure of the tangent bundle. The associated linear Poisson structure on the dual
    (cotangent) bundle provides a framework for the Hamiltonian formulation of the dynamics.
    In this expository  note we show how this correspondence looks like in the affine setting (sections 3, 4).
    First, we introduce the notion of a special affine bundle and its special affine dual (section 2).
    In section 3 we present basic constructions of the geometry of affine values, i.e. the geometry based
    on sections of certain affine bundle, instead of the algebra of smooth functions. In section~4 we
    establish relation between Lie affgebroids (affine version of Lie algebroids) and aff-Jacobi brackets
    on the dual bundle.
    These results are applied in a frame-independent formulation of the time-dependent and
    the Newtonian mechanics (section~5).
%\section{This is the next section}
\section{The affine geometry}
\subsection{Affine spaces and affine bundles}
We start with a short review of basic notions in affine geometry.
We recall some definitions and set the notation. An {\it affine
space} is a triple $(A,V,\za)$, where $A$ is a set, $V$ is a
vector space over a field $\ka$ and $\za$ is a mapping $\za \colon
A \times A\rightarrow V$ such that
      \begin{itemize}
   \item $\za(a_3,a_2) + \za(a_2,a_1) + \za(a_1,a_3) = 0$;
   \item the mapping $\za(\cdot,a) \colon A \rightarrow V$ is bijective for
each $a \in A$.
      \end{itemize}
We shall also write simply $A$ to denote the affine space
$(A,V,\za)$ and $\sV(A)$ to denote its {\it model vector space}
$V$. The {\it dimension of an affine space} $A$ we will call the
dimension of its model vector space $\sV(A)$. A mapping $\zf$ from
the affine space $A_1$ to the affine space $A_2$ is {\it affine}
if there exists a linear map $\zf_\sV: \sV(A_1)\ra\sV(A_2)$ such
that, for every $a\in A_1$ and $u\in \sV(A_1)$,
    $$\zf(a+u)=\zf(a)+\zf_\sV(u).$$
The mapping $\zf_\sV$ is called the {\it linear part} of $\zf$. We
will need also multi-affine, especially bi-affine maps. Let $A$,
$A_1$, $A_2$ be affine spaces. A map
    $$\Phi: A_1\times A_2 \longrightarrow A$$
is called {\it bi-affine} if it is affine with respect to every
argument separately. One can also define linear parts of a
bi-affine map. By $\Phi_{\sV}^1$ we will denote the linear-affine
map
    $$\Phi_{\sV}^1: \sV(A_1)\times A_2\longrightarrow \sV(A)$$
where linear part is taken with respect to the first argument.
Similarly, for the second argument we have the affine-linear map
    $$\Phi_{\sV}^2: A_1\times \sV(A_2)\longrightarrow \sV(A)$$
and finally the {\it bilinear part of\ } $\Phi$:
    $$\Phi_{\sV}:\sV(A_1)\times \sV(A_2)\longrightarrow \sV(A).$$
 For an affine space
$A$ we define its {\it vector dual} $A^\dag$ as a set of all
affine functions  $\phi:A\ra\ka$. The dimension of $A^\dag$ is
greater by $1$ than the dimension of $A$. By $1_A$ we will denote
the element of $A^\dag$ being the constant function on $A$ equal
to $1$. Observe that $A$ is naturally embedded in the space
$\widehat A=(A^\dag)^\ast$. The affine space $A$ can be identified
with the one-codimensional affine subspace of $\widehat A$ of
those linear functions on $A^\dag$ that evaluated on $1_A$ give
$1$. Similarly, $\sV(A)$ can be identified with the vector
subspace of $\widehat A$ of those elements that evaluated on $1_A$
give $0$. The above observation justifies the name {\it vector
hull of $A$} for $\widehat A$.

In the following we will widely use {\it affine bundles} which are
smooth, locally trivial bundles of affine spaces. The notation
will be, in principle, the same for affine spaces and affine
bundles. For instance, $\sV(A)$ denotes the vector bundle which is
the model for an affine bundle $\zz:A\ra M$ over a base manifold
$M$. By $\Sec$ we denote the spaces of sections, e.g. $\Sec(A)$
(or sometimes $\Sec(\zz)$) is the affine space of sections of the
affine bundle $\zz:A\ra M$. The difference of two sections $a$,
$a'$ of the bundle $A$ is a section of $\sV(A)$. Equivalently, we
can add sections of $V(A)$ to sections of $A$:
$$\Sec(A)\ti\Sec(\sV(A))\ni(a,v)\mapsto a+v\in\Sec(A).$$
Every $a\in \Sec(A)$ induces a `linearization'
$$I_a:A\ra\sV(A),\quad a'_p\mapsto a'_p-a(p).$$
The {\it bundle of affine morphisms } $\Aff_M(A_1;A_2)$ from $A_1$
to $A_2$ is a bundle which fibers are spaces of affine maps from
fibers of $A_1$ to fibers of $A_2$ over the same point in $M$. The
space of sections of the bundle we will denote by $\Aff(A_1;A_2)$.
It consists of particular morphisms of affine bundles.
$\Aff_M(A_1;A_2)$ is an affine bundle modelled on
$\Aff_M(A_1;\sV(A_2))$ and $\Aff(A_1;A_2)$ is an affine space
modelled on $\Aff(A_1;\sV(A_2))$. Like for affine maps, we can
define a {\it linear part of an affine morphism}: If $\zf\in
\Aff(A_1;A_2)$ than $\zf_\sV\in \Hom(\sV(A_1);\sV(A_2))$ and for
any $a\in A_1$ and $u\in V(A_1)$, both over the same point in $M$,
we have
    $$\zf_\sv(u)=\zf(a+u)-\zf(a).$$
For an affine bundle $A$ we define also its {\it vector dual}
$A^\dag=\Aff_M(A;M\ti\R)$ and its {\it vector hull}, i.e. the
vector bundle $\wh{A}=(A^\dag)^\ast$.

\subsection{Special affine bundles}
A vector space with a distinguished non-zero vector will be called
{\it a special vector space}. A canonical example of a special
vector space is $(\R,1)$. Another example of a special vector
space is $\bA^\dag=(A^\dag, 1_A)$. A {\it special affine} space is
an affine space modelled on a special vector space. Similarly, a
vector bundle $\bV$ with a distinguished non-vanishing section $v$
will be called {\it a special vector bundle.} An affine bundle
modelled on a special vector bundle will be called {\it a special
affine bundle}. If $\bI$ denotes the special vector space $(\R,1)$
understood also as a special affine space, then $M\times \bI$ is a
special affine bundle over $M$. With some abuse of notation this
bundle will be also denoted by $\bI$. If $\bA_i$, $i=1,2$, denote
special affine bundles modelled on special vector bundles $\bV_i$
with distinguished sections $v_i$ then an affine bundle morphism
$\zf:\bA_1\ra\bA_2$ is a {\it morphism of special affine bundles}
if
$$\zf_{\sV}(v_1)=v_2,$$
i.e. its linear part is a {\it morphism of special vector
bundles}. The set of morphisms of special affine bundles $\bA_i$
will be denoted by $\bAff(\bA_1,\bA_2)$. Similarly, bi-afine
morphisms are called {\it special bi-affine} if they are special
affine with respect to every argument separately.

In the category of special affine bundles there is a canonical
notion of duality. The {\it special affine dual} $\bA^\#$ of a
special vector bundle $\bA=(A,v_A)$ is an affine subbundle in
$A^\dag$ that consists of those affine functions on fibers of $A$
the linear part of which maps $v_A$ to $1$, i.e. those that are
special affine morphisms between $\bA$ and $\bI$:
$$\bA^\#=\bAff(\bA, \bI).$$ The model vector bundle for the dual
$\bA^\#$ is the vector subbundle of $\bA^\dag$ of those functions
which linear part vanishes on $v_A$. We see therefore that $1_A$
is a section of the model vector bundle since its linear part is
identically $0$. Therefore we can consider $\bA^\#$ as a special
affine bundle. The special affine duality is a true duality: in
finite dimensions we have a canonical identification
$$(\bA^\#)^\#\simeq \bA.$$
Given a special affine bundle $\bA=(A,v_A)$, the distinguished
section $v_A$ gives rise to an $\R$ action on $A$, $a_p\mapsto
a_p+tv_A(p)$, whose fundamental vector field will be denoted by
$\chi_\bA$.

\section{The geometry of affine values}
\subsection{The affine phase bundle} In the standard differential
geometry many constructions
are based on the algebra $C^\infty(M)$ of smooth functions on the
manifold $M$. In the AV-geometry we replace $C^\infty(M)$ by the
space of sections of certain affine bundle over $M$. The bundle
$\zz: \Z\ra M$ is a one dimensional affine bundle over the
manifold $M$ modelled on the trivial bundle $M\times \R$. Such a
bundle will be called an {\it AV-bundle. } Since $\Z$ is modelled
on the trivial bundle $M\times \R$, there is an action of $\R$ in
every fiber of $\Z$. Therefore $\Z$ is also a principal
$\R$-bundle. The fundamental vector field induced by the action
will be denoted by $\chi_\Z$. The affine space $\Sec(\Z)$ is
modelled on $C^\infty(M)$. Any section $\sigma$ of $\Z$ induces a
trivialization
    $$ I_\sigma : \Z\ni z\longmapsto  (\zz(z), z-\sigma(\zz(z)))
    \in M\times\R.$$
    The above trivialization induces the identification between
$\Sec(\Z)$ and $C^\infty(M)$. We can go surprisingly far in many
constructions replacing the ring of smooth functions with the
affine space of sections of an AV-bundle. It is possible because
many objects of standard differential geometry (like the cotangent
bundle with its canonical symplectic form) have many properties
that are conserved by certain affine transformations. To build an
AV-analog of the cotangent bundle $\sT^\s M$, let us define an
equivalence relation in the set of all pairs $(m,\zs )$, where $m$
is a point in $M$ and $\sigma $ is a section of $\zz $. Two pairs
$(m,\zs )$ and $(m',\zs ')$ are equivalent if $m' = m$ and
$\xd(\zs ' - \zs)(m) = 0$.  We have identified the section $\zs '
- \zs $ of $pr_M:M\ti\R\ra M$ with a function on $M$ for the
purpose of evaluating the differential $\xd(\zs ' - \zs)(m)$. We
denote by $\sP\Z$ the set of equivalence classes. The class of
$(m,\zs )$ will be denoted by $\bd\sigma (m)$  and will be called
the {\it differential} of $\zs $ at $m$. We will write $\bd$ for
the affine exterior differential to distinguish it from the
standard  $\xd$. We define a mapping ${\sP}\zeta \colon \sP \Z
\rightarrow M$ by $\sP\zz (\bd\zs (m)) = m$. The bundle $\sP\zz$
is canonically an affine bundle modelled on the cotangent bundle
$\zp _M \colon \sT^{\textstyle *} M \rightarrow M$ with the affine
structure
$$\bd\zs _2(m)- \bd\zs _1(m) =
\xd(\zs _2 - \zs _1)(m).
$$
This affine bundle is called the {\it phase bundle} of $\zz$. A
section of $\sP\zz$ will be called an {\it affine 1-form}. Let
$\za\colon M\rightarrow \sP\Z$ be an affine 1-form and let $\zs$
be a section of $\zz$. The differential $\xd(\za -\bd \zs)(m)$
does not depend on the choice of $\zs$ and will be called the {\it
differential of $\za$ at $m$}. We will denote it by $\bd\za(m)$.
The differential of an affine 1-form $\za\in\Sec(\sP\Z)$ is an
ordinary 2-form $\bd\za\in\zW^2(M)$. The corresponding {\it affine
de Rham complex} looks now like
$$\Sec(\Z)\bdra\Sec(\sP\Z)\bdra\zW^2(M)\dra\zW^3(M)\dra\dots$$
and consists of affine maps. This is an {\it affine complex} in
the sense that its linear part is a complex of linear maps.

Like the cotangent bundle $\sT^\ast M$ itself, its AV-analog
$\sP\Z$ is equipped with the canonical symplectic structure. For a
chosen section $\zs$ of $\zz$ we have the isomorphisms
      \bea\nn
   &I_\zs\colon\Z\rightarrow M\times \R,  \\
      &I_{\bd\zs}\colon \sP\Z\rightarrow \sT^{\textstyle *} M,\nn
                                       \eea
   and for two sections $\zs, \zs'$ the mappings $I_{\bd\zs}$ and
   $I_{\bd\zs'}$ differ by the translation by $\xd(\zs -\zs')$, i.e.
      \bea\nn
   &I_{\bd\zs'}\circ I_{\bd\zs}^{-1}\colon \sT^{\textstyle *} M\rightarrow
\sT^{\textstyle *} M  \\
   &\colon\za_m \mapsto \za_m +\xd(\zs -\zs')(m).\nn
                                       \eea
Now we use an affine property of the canonical symplectic form
$\zw_M$ on the cotangent bundle: {\bf translations in $\sT^\s M$
by a closed 1-form are symplectomorphisms}, to conclude that the
two-form $I_{\bd\zs}^{\textstyle *} \zw_M$, where $\zw_M$ is the
canonical symplectic form on $\sT^{\textstyle *} M$, does not
depend on the choice of $\zs$ and therefore it is a canonical
symplectic form on $\sP\Z$. We will denote this form by $\zw_\bZ$.

There is no canonical Liouville 1-form on $\sP\Z$ (in the standard
sense) which is the potential of the canonical symplectic form
$\zw_\Z$ but there is such a form in the affine sense. The details
can be found in \cite{GGU2}.

\subsection{The canonical identification}
Every special affine bundle $\bA=(A, v_{A})$ gives rise to a
certain bundle of affine values $\bAP(\bA)$. The total space of
the bundle $\bAP(\bA)$ is $A$ and the base is $A/\la v_A\ran$ with
the canonical projection from the space to the quotient. The
meaning of the quotient is obvious: the class $[a_p]$ of $a_p\in
A_p$ is the orbit  $\{ a_p+tv_A(p)\}$ of $\chi_\bA$. For the
reason of a further application we choose the distinguished
section $v_{\bAP(\bA)}$ of the model vector bundle characterized
by $\chi_{\bAP(\bA)}=-\chi_{\bA}$. In the space of sections of
$\bAP(\bA)$ one can distinguish {\it affine sections}, i.e. such
sections
$$\bA/\la v_A\ran\ra \bA$$
which are affine morphisms. The space of affine sections will be
denoted by $\Aff\Sec(\bAP(\bA))$.

In the theory of vector bundles there is an obvious identification
between sections of the bundle  $E^\ast$ and functions on $E$
which are linear along fibres. If $\zf$ is a section of $E^\ast$,
then the corresponding function  $\zi_\zf$ is defined by the
canonical pairing
$$\zi_\zf(X)=\langle \zf,X\rangle.$$
In the theory of special affine bundles we have an analog of the
above identification:
$$\Sec(\bA^\#)\simeq\Aff\Sec(\bAP(\bA)),\quad
\bF_\zs\leftrightarrow\zs,$$ where  $\bF_\zs$ is the unique
(affine) function on $\bA$ such that $\chi_\bA(\bF_\zs)=-1$,
$\bF_\zs\circ\zs=0$. In local coordinates $(x^i,s)$ on $A$,
adapted to the structure of the bundle $\bAP(\bA)$, i.e. such that
$(x^i)$ are the coordinates on $A/\la v_A\ran$ and
$\partial_s=-\chi_\bA$ (remember that the fundamental vector field
of $Y$ is the generator of $\exp{(-tY)}$), we have
$$\bF_\zs(x,s)=s-\sigma(x).$$

\section{Brackets on affine bundles}
A {\it Lie affgebra } is an affine space $A$ with a bi-affine
operation
$$[\cdot,\cdot]: A\times A\ra \sV(A)$$
satysfying the following conditions:
  \begin{itemize}
   \item skew-symmetry: $[a_1,a_2]=-[a_2,a_1]$
   \item Jacobi identity:
   $[a_1,[a_2,a_3]]_{\sV}^2+[a_2,[a_3,a_1]]_{\sV}^2+[a_3,[a_1,a_2]]_{\sV}^2=0$.
      \end{itemize}
 In the above formula $[\cdot,\cdot]_{\sV}^2$ denotes the
 affine-linear part of the bracket $[\cdot,\cdot]$, i.e. we take a linear part
 with respect to the second argument. One can show that any
 skew-symmetric bi-affine bracket as above is completely determined by its affine-linear
 part. The detailed description of the concept of the Lie affgebra
 can be found in \cite{GGU1}.

 A {\it Lie affgebroid} is an affine bundle with a Lie affgebra
 bracket on the space of sections
$$[\cdot,\cdot]: \Sec(A)\times \Sec(A)\ra \Sec(\sV(A))$$
and a morphism of affine bundles
$$\rho:A\ra\sT M$$ inducing a map from $\Sec(A)$ into vector
fields on $M$ such that
$$[a,fv]_{\sV}^2=f[a,v]_{\sV}^2+\rho(a)(f)v$$
for any smooth function $f$ on $M$. Note that the same concept has
been introduced by E.~Mart\'\i nez, T.~Mestdag and W.~Sarlet under
the name of {\it affine Lie algebroid} \cite{MMS1}.

\begin{example} Given a fibration $\zx\colon M\ra\R$, take the affine subbundle
$A\subset\sT M$ characterized by  $\zx_*(X)=\pa_t$ for
$X\in\Sec(A)$. Then the standard bracket of vector fields in $A$,
and $\zr\colon A\ra\sT M$ being just the inclusion define on $A$ a
structure of a Lie affgebroid. This is the basic example of the
concept of {\it affine Lie algebroid} developed in
\cite{MMS,MMS1}.
\end{example}

Let us remind that any affine bundle $A$ is canonically embedded
in the vector bundle $\widehat A$ being its vector hull. It turns
out that there is a one-to-one correspondence between Lie
affgebroid structures on $A$ and  Lie algebroid structures on the
vector hull. We have therefore the following theorem (for the
proof see \cite{GGU1}):

\begin{theorem} For an affine bundle $A$ the following are
equivalent:
\begin{itemize}
\item $[\cdot,\cdot]:\Sec(A)\ti\Sec(A)\ra\Sec(\sV(A))$ is a Lie
affgebroid bracket on $A$; \item $[\cdot,\cdot]$ is the
restriction of a Lie algebroid bracket
$[\cdot,\cdot]^\we\colon\Sec(\wh{A})\ti\Sec(\wh{A})\ra\Sec(\wh{A})$.
\end{itemize}
\end{theorem}

Note that $1_A\in \Sec(A^\dag)$ is a closed `one-form' for
$\widehat A$. This is therefore a particular {\it Jacobi
algebroid} as defined in \cite{GM} (or {\it generalized Lie
algebroid} in the terminology of D.~Iglesias and J.~C.~Marrero,
see \cite{IM}).

\begin{definition} Given an AV-bundle $\Z$ over $M$, an {\it aff-Poisson (resp.
aff-Jacobi) bracket} on $\Z$ is a Lie affgebra bracket
$$\{\cdot,\cdot\}\colon\Sec(\Z)\ti\Sec(\Z)\ra C^\infty(M)$$
such that
$$X_\zs=\{\zs,\cdot\}^2_\sv\colon C^\infty(M)\ra
C^\infty(M)$$ is a derivation (resp., a first-order differential
operator) for every $\zs\in\Sec(\Z)$. In the Poisson case $X_\zs$
is a vector field on $M$ called {\it the Hamiltonian vector field
of $\zs$}.
\end{definition}

In the linear case there is a correspondence between Lie algebroid
brackets and linear Poisson structures on the dual bundle. In the
affine setting we have an analog of this correspondence. Let $X$
be a section of $\sV(\bA)$. Since $\sV(\bA)$ is a vector subbundle
in $\wh \bA$, the section X corresponds to a linear function
$\iota^\dag_X$ on $\bA^\dag$. The function $\iota^\dag_X$ is
invariant with respect to the vertical lift of the distinguished
section $1_A$ of $\bA^\dag$. We see that the function
$\iota^\dag_X$ restricted to $\bA^\#$ is constant on fibres of the
projection $\bA^\#\ra\bA^\#/\la 1_A\ran$ and defines an affine
function $\iota^\#_X$ on the base of $\bAP(\bA^\#)$. Hence we have
a canonical identification between
\begin{itemize}
\item  sections $X$ of $\sV(\bA)$, \item linear functions
$\iota_X^\dag$ on $\bA^\dag$ which are invariant with respect to
the vertical lift of $1_A$, \item affine functions $\iota_X^\#$ on
$\bA^\#\slash\langle 1_A\rangle$.
\end{itemize}
Using the above identification we can formulate the following
theorem (for the proof see \cite{GGU2}):

\begin{theorem}
There is a canonical one-to-one correspondence between Lie
affgebroid brackets $[\cdot,\cdot]_\bA$ on $\bA$ and affine
aff-Jacobi brackets $\{\cdot,\cdot\}_{\bA^\#}$ on $\bAP(\bA^\#)$,
uniquely defined by:
$$\{\zs,\zs'\}_{\bA^\#}=\zi^\#_{[\bF_\zs,\bF_{\zs'}]_\bA}.$$
Moreover, $\bA$ is a special Lie affgebroid, i.e. $v$ is a central
section for the Lie algebroid bracket of the vector hull $\wh
\bA$, exactly when the bracket $\{\cdot,\cdot\}_{\bA^\#}$ is
aff-Poisson.
\end{theorem}

\begin{example} Let $\bZ=(Z,v)$ be an $\bAP$-bundle. The tangent
bundle $\sT \bZ$ is equipped with the tangent $\R$-action.
Dividing $\sT\bZ$ by the action we obtain the Atiyah algebroid of
the principal $\R$-bundle $\Z$ which we denote by $\widetilde\sT
\bZ$. It is a special Lie algebroid. The distinguished section of
$\widetilde\sT \bZ$ is represented by the fundamental vector field
$\chi_\bZ$. The AV-bundle $\bAP(\widetilde\sT \bZ)$ is a bundle
over $\sT M$. The special affine dual for $\widetilde\sT \bZ$ is
$\sP \bZ\times \bI$.
%The special affine evaluation between $\sP
%\bZ\times \bI$ and $\widetilde\sT \bZ$ is given by the formula:
%$$\la (\bd\sigma(m), r), \tilde [v]\ran=
%r+([v]-[\sT\sigma(\sT\zz(v))])$$ where the bracket $[\cdot]$
%denotes the equivalence class of $v\in \sT\bZ$ in $\widetilde\sT\bZ$.
The AV-bundle $\bAP((\wt{\sT}\Z)^\#)$ is the trivial bundle
over the affine phase bundle $\sP\Z$ and the correspoding
aff-Poisson bracket is the standard Poisson bracket on $\sP\Z$
associated with the canonical symplectic form $\zw_\Z$ on $\sP\Z$.
\end{example}

\section{Examples} In the following we would like to present two
examples of application of the AV-differential geometry to
classical mechanics. In both examples the Lagrangian mechanics is
associated with a certain special Lie affgebroid $\bA$ and
lagrangians are sections of $\bAP(\bA)$. The Hamiltonian mechanics
in turn is associated with the dual bundle $\bA^\#$ and
hamiltonians are sections of $\bAP(\bA^\#)$.

\subsection{Time-dependent mechanics}

The space of events for the inhomogeneous formulation of the
time-dependent mechanics is the space-time $M$ fibrated over the
time $\T$ being the affine $\R$. The fibration will be denoted by
$\zx$. First-jets of this fibration form the infinitesimal
(dynamical) configuration space. Since there is the distinguished
vector field $\pa_t$ on $\T$, the first-jets of the fibration over
the time can be identified with those vectors tangent to $M$ which
project on $\pa_t$. Such vectors form an affine subbundle $A$ of
the tangent bundle $\sT M$, modelled on the bundle $\sV M$ of
vertical vectors. In the standard formulation the bundle $\sV^\s
M$, dual to the bundle of vectors which are vertical with respect
to the fibration over time, is the phase space for the problem.
The phase space carries a canonical Poisson structure, but
hamiltonian fields for this structure are vertical with respect to
the projection on time, so they cannot describe the dynamics. To
solve this problem one needs the distinguished vector field
$\pa_t$ to be added to the hamiltonian vector field to obtain the
dynamics. This can be done correctly when the fibration over the
time is trivial, i.e. when $M=Q\ti\T$. Then we have the phase
space $\sV^\ast M=\sT^\ast Q\times \T$ and hamiltonian
$$H: \sT^\ast Q\times \R\ni(\alpha,t)\longmapsto H(\alpha,t)\in
\R.$$ The phase dynamics is given by the hamiltonian vector field
$X_{H_t}$ of $H_t=H(\cdot,t)$ calculated separately for every $t$
and corrected by $\partial_t$:
$$X_{H_t}+\partial_t.$$

When the fibration is not trivial one has to choose a reference
vector field that projects onto $\pa_t$. Changing the reference
vector field means changing the hamiltonian.

The Lagrange formalism in the affine formulation originates on the
trivial AV-bundle $\bA=A_0\times\R$ which is a subbundle in $\sT
M\ti \R$. The subbundle $A_0$ is an affine subbundle in $\sT M$ of
those vectors that project onto $\pa_t$. Since $\sT M\times \R$ is
canonically a Lie algebroid with distinguished 1-cocycle $\xd t$,
the bundle $\bA$ is a special Lie affgebroid. Lagrangians are
sections of $\bAP(A)$ and, since the latter bundle is trivial,
they are ordinary functions on $A_0$.

The Hamilton formalism now takes place not on the dual vector
bundle $\sV^\s M$ of $\sV M$, like in the classical approach, but
on the dual AV-bundle $\bAP((\bA)^\#)=\bAP(\bA_0^\dag)$ which can
be recognized as
$$\zz:\sT^\s M\ra\sT^\s M/\la \zx^\ast(dt)\ran $$ and which
carries a canonical aff-Poisson structure induced from the
canonical symplectic Poisson bracket on $\sT^\s M$:
 \be\label{1} \{\zs,\zs'\}\circ\zz=\{\bF_\zs,\bF_{\zs'}\}_M,\ee
where $\{\cdot,\cdot\}_M$ is the canonical Poisson bracket on
$\sT^\s M$. The hamiltonians are sections of the bundle $\zz$.

To compare with the standard approach, let us assume that we have
a trivialization $M=Q\ti\T$ of the space-time into a product of
the space and the time. This induces the decomposition
$$\sT^\s M=\sT^\s Q\ti\sT^\s \T$$
and the bundle $\zz$ becomes
$$\zz: \sT^\s M=\sT^\s Q\ti\sT^\s\T \ra \sT^\s Q\ti\T=\sT^\s M/\la\zx^*(dt)\ran.$$
A time-dependent hamiltonian which is a section of $\zz$ has the
form:
    $$\widehat H(\alpha,t)=(\alpha, t, -H(\alpha,t))$$
and the function on $\sT^\s M$, that corresponds to it, reads
$$ F_{\wh{H}}(\za,t,s)=s+H(\za,t).$$
Finally, we obtain the dynamics
$$\Gamma= X_{\wh{H}}=\zz_*(X_{F_{\wh{H}}})
= \zz_*(X_{H_t}+\pa_t-\frac{\pa H}{\pa t}\pa_s)=X_{H_t}+\pa_t.$$
We see that we have recovered the correct dynamics. However, in
our picture, the term $\pa_t$ is not added `by hand' but it is
generated from $\widehat H$ by means of the aff-Poisson structure.
Of course, if we have no decomposition into space and time, there
is no canonical $\pa_t$ on $M$ and nothing canonical can be added
by hand in the standard approach. This problem disappears in the
aff-Poisson formulation. In this example, hamiltonians are
sections of an AV-bundle and lagrangians are ordinary functions,
however not on a vector but on an affine bundle.

    \begin{remark} The formula (\ref{1}) gives an example of an affine Poisson reduction.
    The standard Poisson reduction deals with a fibration $\zr \colon M\rightarrow N$
    and Poisson brackets $\{\cdot,\cdot\}_M$, $\{\cdot,\cdot\}_N$ on $M$ and $N$, respectively.
    We say that $\{\cdot,\cdot\}_N$ is a reduction of $\{\cdot,\cdot\}_M$ if
    \be\label{2}\{\zr^*f, \zr^*g\}_M = \zr^*\{f,g\}_N.\ee
   In the AV geometry we have two AV-bundles $Z,Y$ with aff-Poisson brackets
   $\{\cdot,\cdot\}_Z$, $\{\cdot,\cdot\}_Y$ on sections of $Z,Y$, respectively.
   A fibration $\zr\colon Z\rightarrow Y$ is assumed to be an AV-morphism.
   The pull-back $\zr^* \zs$ of  a section is well defined and the condition (\ref{2}) can be replaced by
                $$\{\zr^*\zs, \zr^*\zs'\}_Z = \zr^*\{\zs,\zs'\}_Y.$$
    We say that $\{\cdot,\cdot\}_Y$ is an affine reduction (with respect to $\zr$)
    of the aff-Poisson structure $\{\cdot,\cdot\}_Z$.
    In our example (\ref{1}) we have $Z=\sT^\s M \times \R$, $Y= \sT^\s M$, and $\zr(a,r) = a-rdt$.
\end{remark}

\subsection{The inhomogeneous formulation of the Newtonian mechanics}

The Newtonian space-time is a four-dimensional affine space $N$
with the model vector space $\sV(N)$. There is a distinguished
non-zero covector $\tau\in\sV(N)^\ast$ that measures time
intervals between events:
$$\Delta t(x_1,x_2)=\la\tau, x_2-x_1\ran$$
There is also an Euclidean metric $g$ defined on the subspace
$E_0=\ker \tau$. It serves for measuring the spatial distance
between symultaneous events. In the standard approach the dynamics
of a massive particle is described in an inertial frame. An
inertial frame is a class of inertial observers that move in
space-time with the same constant velocity. A level-$1$ set of
$\tau$ represents velocities of inertial observers and velocities
of particles.

When we fix an inerial frame $u$, the model vector space $\sV(N)$
splits into $E_0\times \R$:
$$\sV(N) \ni v\longmapsto (v-\la\tau,v\ran u, \,\,\,\la\tau,v\ran)$$
and the dual $\sV(N)^\ast$ splits into $E_0^\ast\times\R$:
$$\sV(N)^\ast\ni \alpha \longmapsto(\iota^\ast(\alpha),\,\,\, \la\alpha, u\ran)\in
E_0^\ast\ti \R,$$ where $\iota^\ast$ is a canonical projection
which is dual to the embedding $\iota: E_0\hookrightarrow\sV(N)$.
In the standard approach the phase space of a massive particle is
$N\times E_0^\ast$ that is interpreted as the bundle which is dual
to the bundle $N\times E_0$ of vectors tangent to $N$ and vertical
with respect to the projection on the time given by $\tau$. The
phase of a particle with a mass $m$ moving with a relative
velocity $v\in E_0$ in the field of forces with a potential
$\zf(x)$ is $(x,p)=(x,mg(v))\in N\ti E_0^*$ and the Hamilton
function is
$$H_u=H_u(x,p)=\frac{1}{2m}p^2+\zf(x),$$
with $p^2=\la p,g^{-1}(p)\ran$.  The dynamics is generated from
the hamiltonian $H_u$ by means of the canonical Poisson structure
of $N\times E_0^\ast$ obtained by reduction from $\sT^\ast N\simeq
N\times \sV(N)^\ast$. We meet here the same problem as in the
previous example: the Hamiltonian vector field generated from
$H_u$ is vertical with respect to the projection on the time. Here
we have to add the constant vector field equal to $u$ at every
point to obtain the correct phase equations. Using the AV-geometry
we can get rid of the inertial observer and find the correct
aff-Poisson structure that generates the correct phase equations.

Analyzing transformation rules for the energy and the momentum
while changing the observer, we obtain the gauge-independent
(observer-independent) phase AV-bundle
$$\widehat\sP\longrightarrow\sP$$ where
$$\wh{\sP}=(N\ti E_0^*\ti\R\ti E_1)/E_0\quad{\rm and} \quad\sP=(N\ti E_0^*\ti
E_1)/E_0$$ with respect to the gauge equivalence
\beas(x,p,s,u'+v)&\sim&(x',p',s,u')\Leftrightarrow x=x', \\
&& p=p'+mg(v), \\ && s=s'+\la p,v\ran+\frac{1}{2}m\la
g(v),v\ran.\eeas The phase AV-bundle comes together with a
canonical aff-Poisson structure. This structure is inherited from
the AV bundle $\sT^\ast N\ra \sT^\ast N/\la\tau\ran\simeq N\times
E_0^\ast$. We have seen in the previous example that such a bundle
has a canonical aff-Poisson structure coming from the symplectic
structure of the cotangent bundle. What we claim here is that the
aff-Poisson structure survives the procedure of collecting phase
bundles for all inertial frames and then dividing by
transformation rules.

The potential $\zf(x)$ gives rise to a well-defined Hamiltonian
section $\wh{H}$ of $\bAP(\wh{\sP})$ which for the trivialization
given by the inertial frame $u$
$$I_{u}\colon\,\,\wh{\sP}\longrightarrow N\ti E_0^*\ti\R,$$
reads $$H(x,p)=\frac{1}{2m}p^2+\zf(x).$$ The observed dynamics on
the phase space $\sP\simeq N\times E_0$ is
\beas\dot{x}&=&\left(g^{-1}\left(\frac{p}{m}\right)\right)+u,\\
\dot{p}&=&-\frac{\pa\zf}{\pa q}(q,t).\eeas If we fix additionally
an inertial observer from the class characterized by $u$,
specifying a point $x_0$ in $N$, then $N$ splits into
$E_0\times\R$ and the dynamics on $E_0\times\R\times E_0^\ast$,
with elements denoted by $(q,t,p)$, reads
\beas \dot{q}&=&\left(g^{-1}\left(\frac{p}{m}\right)\right),\\
\dot{t}&=& 1,\\
\dot{p}&=&-\frac{\pa\zf}{\pa q}(q,t),\eeas that agrees with the
Newtonian picture.

%%%%%%%%%%%%%%%%%%%%%%%%%%%%%%%%%%%%%
%%%%%%%%%%%%%%%%%%%%%%%%%%%%%%%%%%%%%%%%%%%%%%%%%%%%%%%%%%%%%%%%%%%%%%%
%%%%%%%%%%%%%%%%%%%%%%%%%%%%%%%%%%%%%%%%%%%%%%%%%%%%%%%%%%%%%%%%%%%%
%%%%%          The bibliography
%%%%%%%%%%%%%%%%%%%%%%%%%%%%%%%%%%%%%%%%%%%%%%%%%%%%%%%%%%%%%%%%%%%%%

%%%%%%%%%%%%%%%%%%%%%%%%%%%%%%%%%%%%%%%%%%%%%%%%%%%%%%%%%%%%
%%%%%%%%%%%%%%%%%%%%%%%%%%%%%%%%%%%%%%%%%%%
%% The address
\noindent
Katarzyna  Grabowska \\
        Division of Mathematical Methods in Physics\\
                University of Warsaw \\
                Ho\.za 69, 00-681 Warszawa, Poland \\
                 {\tt konieczn@fuw.edu.pl} \\
                        \\
        Janusz Grabowski \\
        Institute of Mathematics \\
                Polish Academy of Sciences \\
                ul. \'Sniadeckich 8, P.O.Box 21, 00-956 Warszawa 10, Poland\\
                {\tt jagrab@impan.gov.pl} \\
                        \\
               Pawe\l\ Urba\'nski \\
                Division of Mathematical Methods in Physics \\
                University of Warsaw \\
                Ho\.za 69, 00-681 Warszawa, Poland \\
                {\tt urbanski@fuw.edu.pl}
%%%%%%%%%%%%%%%%%%%%%%%%%%%%%%%%%%%%%%%%%%%%%%%%%%%%%%%%%%%%%%%%%

\label{lastpage}
\end{document}